\documentclass[10pt]{amsart}
\usepackage{amsfonts}
\usepackage{ifthen}
\usepackage{amsthm}
\usepackage{amsmath}
\usepackage{graphicx}
\usepackage{amscd,amssymb,amsthm}
\usepackage{enumerate}
\usepackage{epstopdf}
\usepackage{mathrsfs}
\usepackage{hyperref}

\newcounter{minutes}
\divide\time by 60
\newcounter{hours}
\multiply\time by 60 \addtocounter{minutes}{-\time}

\title{Bounds for the asymptotic order parameter of the stochastic Kuramoto model}

\author[I. Mez\H o]{Istv\'an Mez\H o${\dag}$}
\address{Department of Mathematics, Nanjing University of Information Science and
Technology, 5 Panxin Rd, Pukou, Nanjing, Jiangsu, P.R. China}
\email{istvanmezo81@gmail.com}

\author[\'A. Baricz]{\'Arp\'ad Baricz${\ddag}$}
\address{Institute of Applied Mathematics, \'Obuda University, 1034 Budapest, Hungary}
\address{Department of Economics,  Babe\c{s}-Bolyai University, Cluj-Napoca 400591, Romania}
\email{bariczocsi@yahoo.com}

\thanks{${\dag}$The research of Istv\'an Mez\H{o} was supported by the Scientific Research Foundation of Nanjing University of Information Science \& Technology, the Startup Foundation for Introducing Talent of NUIST, Project no. S8113062001, and the National Natural Science Foundation for China, Grant no. 11501299.}

\thanks{${\ddag}$The research of \'Arp\'ad Baricz was supported by the J\'anos Bolyai Research Scholarship of the Hungarian Academy of Sciences. The work of \'A. Baricz was initiated during his visit in October 2015 to Department of Mathematics of Nanjing University of Information Science and Technology, to which this author is grateful for hospitality.}

\begin{document}

\def\thefootnote{}
\footnotetext{ \texttt{File:~\jobname .tex,
          printed: \number\year-0\number\month-\number\day,
          \thehours.\ifnum\theminutes<10{0}\fi\theminutes}
} \makeatletter\def\thefootnote{\@arabic\c@footnote}\makeatother

\keywords{Stochastic Kuramoto model, asymptotic order parameter, modified Bessel functions, Tur\'an type inequalities, approximation, Lagrange inversion, monotonicity properties, bounds.}
\subjclass[2010]{39B62, 33C10.}

\maketitle

\begin{center}
Dedicated to Bor\'oka, Eszter and Kopp\'any
\end{center}

\begin{abstract}
Tur\'an type inequalities for modified Bessel functions of the first kind are used to deduce some
sharp lower and upper bounds for the asymptotic order parameter of the stochastic Kuramoto model. Moreover, approximation from the
Lagrange inversion theorem and a rational approximation are given for the asymptotic order parameter.
\end{abstract}

\section{Introduction}

The Kuramoto model describes the phenomenon of collective synchronization, more precisely it
describes how the phases of coupled oscillators evolve in time, see \cite{kur} for more details.
Recently Bertini, Giacomin and Pakdaman \cite{bertini} were able to review some results on the
Kuramoto model from a statistical mechanics standpoint and they gave in particular necessary and
sufficient conditions for reversibility. In order to do this Bertini, Giacomin and Pakdaman \cite[p. 278]{bertini}
deduced some lower and upper bounds for the asympotic order parameter, which involves the
modified Bessel functions of the first kind of order zero and one. A few years later Sonnenschein and
Schimansky-Geier \cite{bernard} obtained the asymptotic order parameter in closed form, which suggested
a tighter upper bound for the corresponding scaling. Moreover, they elaborated the Gaussian
approximation in complex networks with distributed degrees. In their study Sonnenschein and
Schimansky-Geier \cite[p. 3]{bernard} proposed another upper bound for the asymptotic order parameter,
but they presented their result without mathematical proof. All the same, by using Bernoulli's inequality
they verified that their upper bound is better than the upper bound of Bertini, Giacomin and Pakdaman \cite{bertini}.
In this paper our aim is to make a contribution to this subject by showing the followings:

\begin{enumerate}
\item[$\bullet$] The bounds presented in the above mentioned papers are correct and their proofs are based on some Tur\'an type inequalities for modified Bessel functions of the first kind.
\item[$\bullet$] The constants in the upper bounds presented by Bertini, Giacomin, Pakdaman \cite{bertini} and Sonnenschein, Schimansky-Geier \cite{bernard} are the best, and thus their bounds cannot be improved.
\item[$\bullet$] The results presented in the above mentioned papers can be extended to modified Bessel functions of the first kind of arbitrary order, based on some interesting new and recently discovered Tur\'an type inequalities for modified Bessel functions of the first kind.
\item[$\bullet$] It is possible to obtain another approximation for the asymptotic order parameter (than in the above mentioned papers) by means of the Lagrange's inversion theorem and also a rational approximation.
\end{enumerate}

As far as we know the above mentioned subject was not studied yet in details from the mathematical point of view and we believe that the obtained results may be useful for the people working in statistical physics.

\section{Bounds for the asymptotic order parameter}
\setcounter{equation}{0}

In this section our aim is to discuss, complement and extend the results from \cite{bertini,bernard} concerning bounds for the asymptotic order parameter of the stochastic Kuramoto model. Some new and recently discovered Tur\'an type inequalities for modified Bessel functions of the first kind play an important role in this section. For more details on Tur\'an type inequalities for modified Bessel functions of the first kind we refer to \cite{baricz} and to the references therein.

\subsection{An alternative proof of a result on asymptotic order parameter.} Let us consider the transcendental equation $r=\Psi(2Kr),$ where $K>1,$ $\Psi(x)=I_1(x)/I_0(x)$ and $I_1,$ $I_0$ stand for the modified Bessel function of the first kind of order $1,$ and $0,$ respectively. Recently, Bertini, Giacomin and Pakdaman \cite{bertini} in order to prove their main result about the spectrum of a self-adjoint linear operator, presented the inequalities
\begin{equation}\label{ineq1}\sqrt{1-\frac{1}{K}}<r<\sqrt{1-\frac{1}{2K}}.\end{equation}
The clever proof of the left-hand side of \eqref{ineq1} was based on the well-known Tur\'an type inequality
$$I_1^2(x)-I_0(x)I_2(x)>0.$$ In what follows we would like to show that in fact the right-hand side of \eqref{ineq1} is also equivalent to a Tur\'an type inequality involving modified Bessel functions of the first kind. To proceed, we use the same notation as in \cite{bertini}. To prove the right-hand side of \eqref{ineq1} we need to show that
$$r^2+\frac{1}{2K}-1=\Psi^2(2Kr)+\frac{\Psi(2Kr)}{2Kr}-1<0,$$
that is, for $x>0$ we have
\begin{equation}\label{ineqpsi1}\Psi^2(x)+\frac{1}{x}\Psi(x)-1<0.\end{equation}
Now, by applying the identity \cite[eq. 2.6]{bertini}
\begin{equation}\label{quot1}\frac{I_1(x)}{I_0(x)}=\frac{x}{2}\left(1+\frac{x}{2}\frac{I_2(x)}{I_1(x)}\right)^{-1}\end{equation}
we obtain $$1-\frac{1}{x}\Psi(x)=\left({1+\frac{1}{x}\frac{I_1(x)}{I_2(x)}}\right)\left({1+\frac{2}{x}\frac{I_1(x)}{I_2(x)}}\right)^{-1},$$
which implies that \eqref{ineqpsi1} is equivalent to
$$\Psi^2(x)\left({1+\frac{2}{x}\frac{I_1(x)}{I_2(x)}}\right)<{1+\frac{1}{x}\frac{I_1(x)}{I_2(x)}},$$
which by means of the recurrence relation \begin{equation}\label{rec}xI_0(x)-xI_2(x)=2I_1(x),\end{equation} is equivalent to the Tur\'an type inequality
$$I_1^2(x)-I_0(x)I_2(x)<\frac{1}{x}I_0(x)I_1(x).$$
But, in view of the well-known Soni inequality $I_1(x)<I_0(x),$ the above Tur\'an type inequality is a consequence of the stronger inequality \cite[eq. 2.5]{baricz}
$$I_1^2(x)-I_0(x)I_2(x)<\frac{1}{x}I_1^2(x).$$
Since all of the above inequalities are valid for $x>0$ it follows that the right-hand side of \eqref{ineq1} is valid.

\subsection{The proof of a claimed result on asymptotic order parameter.}
Recently, Sonnenschein and Schimansky-Geier \cite{bernard} proposed (without proof) an improvement of the right-hand side of \eqref{ineq1} as follows
\begin{equation}\label{ineq2}r<\sqrt[4]{1-\frac{1}{K}}.\end{equation}
In the sequel we present a proof of \eqref{ineq2}, which is based also on a Tur\'an type inequality. Note that to prove \eqref{ineq2} we need to show that
$$r^4+\frac{1}{K}-1=\Psi^4(2Kr)+\frac{2\Psi(2Kr)}{2Kr}-1<0,$$
that is, for $x>0$ we have
\begin{equation}\label{ineqpsi2}\Psi^4(x)+\frac{2}{x}\Psi(x)-1<0.\end{equation}
Now, by applying the identity \eqref{quot1} we obtain $$1-\frac{2}{x}\Psi(x)=\left({1+\frac{2}{x}\frac{I_1(x)}{I_2(x)}}\right)^{-1},$$
which implies that \eqref{ineqpsi2} is equivalent to
$$\Psi^4(x)\left({1+\frac{2}{x}\frac{I_1(x)}{I_2(x)}}\right)<1,$$
which by means of the recurrence relation \eqref{rec} is equivalent to the Tur\'an type inequality
\begin{equation}\label{edin}I_1^4(x)<I_0^3(x)I_2(x).\end{equation}
Since all of the above inequalities are valid for $x>0$ it follows that indeed the right-hand side of \eqref{ineq1} is valid.
The Tur\'an type inequality \eqref{edin} is the limiting case of the next Tur\'an type inequality (see \cite[p. 592]{baed}) when $\nu\to-1$
$$I_{\nu+1}^3(x)I_{\nu+3}(x)>I_{\nu}(x)I_{\nu+2}^3(x),\ \ x>0,\ \nu>-1,$$
and it was shown by Idier and Collewet \cite{idier} that we can take the limit and the inequality \eqref{edin} remains true.

\subsection{Sharpness of the results concerning the asymptotic order parameter} It is important to mention here that the powers in the left-hand side of the inequality \eqref{ineq1}, and in \eqref{ineq2}, that is
\begin{equation}\label{ineq3}\sqrt{1-\frac{1}{K}}<r<\sqrt[4]{1-\frac{1}{K}},\end{equation}
are the best possible. To show this observe that the inequality \eqref{ineq3} is equivalent to
$$2<\frac{\log\left(1-\frac{1}{K}\right)}{\log r}<4\ \ \ \ \mbox{or to}\ \ \ \ 2<\frac{\log\left(1-\frac{2}{x}\Psi(x)\right)}{\log\Psi(x)}<4.$$
Here we used the inequality $I_1(x)<I_0(x),$ where $x>0,$ which shows that $\Psi$ maps $(0,\infty)$ into $(0,1),$ and thus $\log\Psi(x)<0$ for $x>0.$ Now we show that in the above inequalities the constants $2$ and $4$ are the best possible in the sense that the inequality
$$\alpha<\frac{\log\left(1-\frac{2}{x}\Psi(x)\right)}{\log\Psi(x)}<\beta$$
is valid for all $x>0$ with the optimal parameters $\alpha=2$ and $\beta=4.$ This implies that in \eqref{ineq3} the powers $\frac{1}{2}$ and $\frac{1}{4}$ are indeed the best possible. Thus, let
\[\lambda(x)=\frac{\log\left(1-\frac2x\Psi(x)\right)}{\log\Psi(x)}.\]
We prove that $\alpha=\lim\limits_{x\to0}\lambda(x)=2.$ Because of the presence of the logarithm, $\lambda(x)$ has no Taylor series around $x=0$, but it still can be expanded as follows
\[\lambda(x)=\frac{\log (8)-2 \log (x)}{\log (2)-\log (x)}+x^2\frac{2 \log (x)+4 \log (2)-3 \log (8)}{24 (\log (2)-\log (x))^2}+\cdots,\]
from where the above limit immediately comes, as we stated. Note that by using the Mittag-Leffler expansion of $\Psi(x)$ it is also possible to obtain the above limit. Namely, since
$$\Psi(x)=\frac{I_1(x)}{I_0(x)}=\sum_{n\geq1}\frac{2x}{x^2+j_{0,n}^2},$$
where $j_{0,n}$ stands for the $n$th positive zero of the Bessel function $J_0,$ by using the Bernoulli-l'Hospital's rule twice we obtain that
$$\lim_{x\to0}\lambda(x)=8\lim_{x\to0}\frac{\displaystyle\sum_{n\geq1}\frac{j_{0,n}^2-x^2}{(x^2+j_{0,n}^2)^2}\displaystyle\sum_{n\geq1}\frac{x}{(x^2+j_{0,n}^2)^2}+
\displaystyle\sum_{n\geq1}\frac{x}{x^2+j_{0,n}^2}\displaystyle\sum_{n\geq1}\frac{j_{0,n}^2-3x^2}{(x^2+j_{0,n}^2)^3}}{4\displaystyle\sum_{n\geq1}\frac{2x}{(x^2+j_{0,n}^2)^2}
\displaystyle\sum_{n\geq1}\frac{j_{0,n}^2-x^2}{(x^2+j_{0,n}^2)^2}+\left(1-4\displaystyle\sum_{n\geq1}\frac{1}{x^2+j_{0,n}^2}\right)\displaystyle\sum_{n\geq1}\frac{2x(x^2-3j_{0,n}^2)}{(x^2+j_{0,n}^2)^3}}=2.$$

Now, we are going to prove that the best constant $\beta$ equals to $\beta=\lim\limits_{x\to\infty}\lambda(x)=4.$ The well known asymptotic estimation
\begin{equation}\label{assym}I_\nu(x)=\frac{e^x}{\sqrt{2\pi x}}\left(1-\frac{4\nu^2-1}{8x}+O\left(\frac{1}{x^2}\right)\right)\end{equation}
yields that as $x$ grows
\[\Psi(x)=\frac{8x-4+1+O\left(\frac{1}{x}\right)}{8x+1+O\left(\frac{1}{x}\right)}.\]
Substituting this into the definition of $\lambda(x)$ we can see that asymptotically it equals to
\[\lambda(x)=\frac{\log(x-2+O\left(\frac{1}{x}\right))-\log(x)}{\log\left(8x-4+1+O\left(\frac{1}{x}\right)\right)-\log\left(8x+1+O\left(\frac{1}{x}\right)\right)}.\]
The differences of logarithms both in the nominator and denominator can be expanded at infinity by the expansion
\[\log(ax+b)-\log(cx+d)=\log(a)-\log(c)+\left(\frac{b}{a}-\frac{d}{c}\right)\frac1x+O\left(\frac{1}{x^2}\right),\quad ac\neq0.\]
What we get is the following
\[\lambda(x)=\frac{\frac2x+O\left(\frac{1}{x^2}\right)}{\left(\frac{1}{2}\right)\frac1x+O\left(\frac{1}{x^2}\right)}=
\frac{4+O\left(\frac{1}{x}\right)}{1+O\left(\frac{1}{x}\right)}\]
for large $x$, from which the result follows.

\subsection{Extension of \eqref{ineq1} to the general case} Let us consider the function $\Psi_{\nu}:(0,\infty)\to(0,\infty),$ defined by $\Psi_{\nu}(x)=I_{\nu+1}(x)/I_{\nu}(x),$ where $I_{\nu}$ stands for the modified Bessel function of the first kind of order $\nu.$ We are going to show that it is possible to extend the right-hand side of \eqref{ineq1} to the general case. Thus, we consider the transcendental equation $r=\Psi_{\nu}(2Kr)$ and we show that the next inequality holds true for all $K>\nu+1,$ $\nu\geq 0$ and $r>0$
\begin{equation}\label{ineq4}r<\sqrt{1-\frac{1}{2K}}.\end{equation}
However, we first show that if $\nu\ge 0$ and $K>\nu+1,$ then the equation
\begin{equation}
r-\Psi_\nu(2Kr)=0\label{Kcond}
\end{equation}
has a positive real solution. Observe that our equation has a trivial solution at $r=0$. Moreover, the function $r-\Psi_\nu(2Kr)$ tends to infinity as $r$ grows. Hence, if the tangent of the continuous function $r-\Psi_\nu(2Kr)$ is negative at the origin then \eqref{Kcond} surely has a positive real solution. The derivative is as follows
\[(r-\Psi_\nu(2Kr))'=1+\frac{K I_{\nu+1}(2 K r) (I_{\nu-1}(2 K r)+I_{\nu+1}(2 K r))}{I_{\nu}(2 K r){}^2}-\frac{K (I_{\nu}(2 K r)+I_{\nu+2}(2 K r))}{I_{\nu}(2 K r)}.\]
We need to take the limit when $r=0$. This can be done by using the Bernoulli-l'Hospital's rule. It comes after some simplification that
\[\left.(r-\Psi_\nu(2Kr))'\right|_{r=0}=1-\frac{K}{\nu+1},\]
from where it follows that in the case when $1-\frac{K}{\nu+1}<0$ then \eqref{Kcond} has a positive real solution.

Now, to prove the inequality \eqref{ineq4} we need to show that
$$r^2+\frac{1}{2K}-1=\Psi_{\nu}^2(2Kr)+\frac{\Psi_{\nu}(2Kr)}{2Kr}-1<0,$$
that is, for $x>0$ we have
\begin{equation}\label{ineqpsigen}\Psi_{\nu}^2(x)+\frac{1}{x}\Psi_{\nu}(x)-1<0.\end{equation}
By applying the identity
\begin{equation}\label{quotgen}\frac{I_{\nu+1}(x)}{I_{\nu}(x)}=\frac{x}{2}\left(\nu+1+\frac{x}{2}\frac{I_{\nu+2}(x)}{I_{\nu+1}(x)}\right)^{-1}\end{equation}
we obtain $$1-\frac{1}{x}\Psi_{\nu}(x)=\left({1+\frac{2\nu+1}{x}\frac{I_{\nu+1}(x)}{I_{\nu+2}(x)}}\right)
\left({1+\frac{2(\nu+1)}{x}\frac{I_{\nu+1}(x)}{I_{\nu+2}(x)}}\right)^{-1},$$
which implies that \eqref{ineqpsigen} is equivalent to
$$\Psi_{\nu}^2(x)\left({1+\frac{2(\nu+1)}{x}\frac{I_{\nu+1}(x)}{I_{\nu+2}(x)}}\right)<{1+\frac{2\nu+1}{x}\frac{I_{\nu+1}(x)}{I_{\nu+2}(x)}},$$
which by means of the recurrence relation \begin{equation}\label{recgen}xI_{\nu}(x)-xI_{\nu+2}(x)=2(\nu+1)I_{\nu+1}(x),\end{equation} is equivalent to the Tur\'an type inequality
$$I_{\nu+1}^2(x)-I_{\nu}(x)I_{\nu+2}(x)<\frac{2\nu+1}{x}I_{\nu}(x)I_{\nu+1}(x).$$
But, in view of the well-known Soni inequality $I_{\nu+1}(x)<I_{\nu}(x),$ $\nu>-\frac{1}{2},$ $x>0,$ the above Tur\'an type inequality is a consequence of the stronger inequality \cite[eq. 2.5]{baricz}
\begin{equation}\label{turanb}I_{\nu+1}^2(x)-I_{\nu}(x)I_{\nu+2}(x)<\frac{1}{x}I_{\nu+1}^2(x),\end{equation}
which holds for $\nu\geq-\frac{1}{2}$ and $x>0.$ Since all of the above inequalities are valid for $x>0$ and $\nu\geq0,$ it follows that indeed the inequality \eqref{ineq4} is valid.

Moreover, it can be shown that the left-hand side of \eqref{ineq1} can be also extended to the general case when $\nu\geq\frac{1}{2},$ and the resulting inequality is reversed (comparative to the left-hand side of \eqref{ineq1}) and improves the inequality \eqref{ineq4}. To show the inequality
\begin{equation}\label{ineq5}r<\sqrt{1-\frac{1}{K}}\end{equation}
it is enough to show that for $x>0$ we have
\begin{equation}\label{ineqpsigen2}\Psi_{\nu}^2(x)+\frac{2}{x}\Psi_{\nu}(x)-1<0.\end{equation}
By using the above steps in view of \eqref{quotgen} we obtain that
$$1-\frac{2}{x}\Psi_{\nu}(x)=\left({1+\frac{2\nu}{x}\frac{I_{\nu+1}(x)}{I_{\nu+2}(x)}}\right)
\left({1+\frac{2(\nu+1)}{x}\frac{I_{\nu+1}(x)}{I_{\nu+2}(x)}}\right)^{-1},$$
which implies that \eqref{ineqpsigen2} is equivalent to
$$\Psi_{\nu}^2(x)\left({1+\frac{2(\nu+1)}{x}\frac{I_{\nu+1}(x)}{I_{\nu+2}(x)}}\right)<{1+\frac{2\nu}{x}\frac{I_{\nu+1}(x)}{I_{\nu+2}(x)}},$$
which by means of the recurrence relation \eqref{recgen} is equivalent to the Tur\'an type inequality
\begin{equation}\label{turaninter}I_{\nu+1}^2(x)-I_{\nu}(x)I_{\nu+2}(x)<\frac{2\nu}{x}I_{\nu}(x)I_{\nu+1}(x).\end{equation}
But, in view of the above Soni inequality the above Tur\'an type inequality is a consequence of the stronger inequality \eqref{turanb} when $\nu\geq\frac{1}{2}.$

Numerical experiments suggest that when $\nu\in\left(0,\frac{1}{2}\right)$ the equation
$$I_{\nu+1}^2(x)-I_{\nu}(x)I_{\nu+2}(x)-\frac{2\nu}{x}I_{\nu}(x)I_{\nu+1}(x)=0$$
has a solution which depends on $\nu.$ Since \eqref{ineq5} is equivalent to \eqref{turaninter}, this
implies that \eqref{ineq5} is not true for all $K>\nu+1,$ $\nu\in\left(0,\frac{1}{2}\right)$ and $r>0.$

\subsection{Sharpness of the extension \eqref{ineq5}} Let us consider the extension of $\lambda,$ that is,
\[\lambda_\nu(x)=\frac{\log\left(1-\frac2x\Psi_\nu(x)\right)}{\log\Psi_\nu(x)}.\]
We are going to prove that the best constant $\beta_\nu$ for which $\lambda_{\nu}(x)<\beta_{\nu}$ for $x>0$ and $\nu\geq0,$ equals to
\[\beta_\nu=\lim_{x\to\infty}\lambda_\nu(x)=\frac{4}{2\nu+1}.\]
In particular, $\beta_0=\beta=4.$ The steps are the same as in the special case above when $\nu=0.$ The asymptotic estimation \eqref{assym} yields that as $x$ grows
\begin{equation}\label{assymp2}\Psi_\nu(x)=\frac{8x-4(\nu+1)^2+1+O\left(\frac{1}{x}\right)}{8x-4\nu^2+1+O\left(\frac{1}{x}\right)}.\end{equation}
Substituting this into the definition of $\lambda_\nu(x)$ we can see that asymptotically it equals to
\[\lambda_\nu(x)=\frac{\log(x-2+O\left(\frac{1}{x}\right))-\log(x)}
{\log\left(8x-4(\nu+1)^2+1+O\left(\frac{1}{x}\right)\right)-\log\left(8x-4\nu^2+1+O\left(\frac{1}{x}\right)\right)},\]
and consequently
\[\lambda_\nu(x)=\frac{\frac2x+O\left(\frac{1}{x^2}\right)}{\left(\frac{(\nu+1)^2}{2}-\frac{\nu^2}{2}\right)\frac1x+O\left(\frac{1}{x^2}\right)}=\frac{4+O\left(\frac{1}{x}\right)}{2\nu+1+O\left(\frac{1}{x}\right)}\]
for large $x$, from which the result follows. The above discussion actually shows that when $\nu\geq \frac{1}{2}$ the extension \eqref{ineq5} is far from being the best possible one. The next natural extension of \eqref{ineq3} improves the extension \eqref{ineq5} when $\nu\geq \frac{1}{2}$ and the power $\frac{\nu}{2}+\frac{1}{4}$ appearing in the inequality is the best possible
\begin{equation}\label{ineq6}r<\sqrt[4]{\left(1-\frac{1}{K}\right)^{2\nu+1}}.\end{equation}

\begin{center}
\begin{figure}[!ht]
   \centering
       \includegraphics[width=13cm]{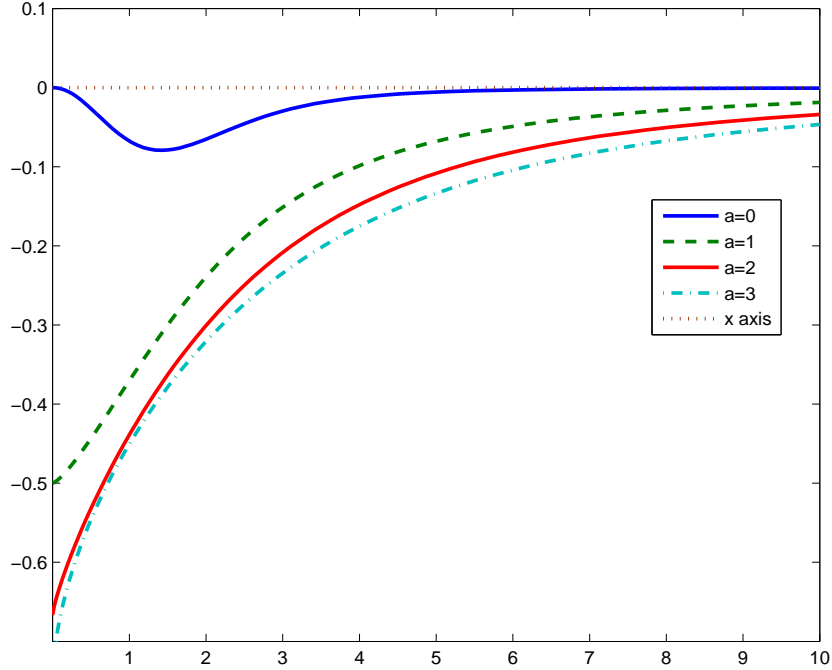}
       \caption{The graph of the function $x\mapsto \sqrt[2a+1]{\Psi_{a}^4(x)}+\frac{2}{x}\Psi_{a}(x)-1$ on $[0,10]$ in the case when $a\in\{0,1,2,3\}.$}
       \label{fig1}
\end{figure}
\end{center}

\subsection{Sharp extension of \eqref{ineq1} to the general case} Now, we are going to show that \eqref{ineq6} is valid for all $K>\nu+1,$ $\nu\geq 0.3$ and $r>0.$ Note that the inequality \eqref{ineq6} is equivalent to
\begin{equation}\label{ineq7}\sqrt[2\nu+1]{\Psi_{\nu}^4(x)}+\frac{2}{x}\Psi_{\nu}(x)-1<0,\end{equation}
where $\nu\geq 0,$ $x>0.$ Using the Amos type bound $\Psi_{\nu}(x)<\Omega_{\nu}(x),$ where $\nu\geq0,$ $x>0$ and \cite[p. 94]{hg}
$$\Omega_{\nu}(x)=\frac{x}{\sqrt{x^2+\left(\nu+\frac{1}{2}\right)\left(\nu+\frac{3}{2}\right)}+\nu+\frac{1}{2}},$$
we obtain that
$$\sqrt[2\nu+1]{\Psi_{\nu}^4(x)}+\frac{2}{x}\Psi_{\nu}(x)-1<\sqrt[2\nu+1]{\Omega_{\nu}^4(x)}+\frac{2}{x}\Omega_{\nu}(x)-1<0,$$
where $\nu\geq0.3$ and $x>0.$ Here we used the fact that, based on numerical experiments, the smallest value of $\nu$ for which the expression $\sqrt[2\nu+1]{\Omega_{\nu}^4(x)}+\frac{2}{x}\Omega_{\nu}(x)-1$ is still negative for each $x>0$ is $0.3.$ When $\nu=0.3$ the above expression tends to zero as $x\to0.$ We believe, but were unable to prove that \eqref{ineq6} is also valid when $\nu\in(0,0.3),$ $K>\nu+1$ and $r>0.$ Numerical experiments (see also Fig. \ref{fig1}) strongly suggest the validity of the above claim.

It is also worth to mention that the inequality \eqref{ineq7} is actually equivalent to the Tur\'an type inequality
$$\sqrt[2\nu+1]{\left(\frac{I_{\nu+1}(x)}{I_{\nu}(x)}\right)^{4}}\cdot\frac{I_{\nu}(x)}{I_{\nu+2}(x)}<1+\frac{2\nu}{x}\frac{I_{\nu+1}(x)}{I_{\nu+2}(x)},$$
where $\nu\geq 0.3$ and $x>0.$ This Tur\'an type inequality is new, and we believe, but were unable to prove that it is true also when $\nu\in(0,0.3)$ and $x>0.$ The case $\nu=0$ has been already considered in \eqref{edin}.

\subsection{Another sharp extension of \eqref{ineq1} to the general case} In this subsection our aim is to propose that it would be possible to extend the inequality \eqref{ineq1} in another way such that to keep the sharpness. We believe that if $\nu\geq0,$ $K>\nu+1$ and $r>0,$ then we have
\begin{equation}\label{ineq8}r<\sqrt[4(\nu+1)]{\left(1-\frac{\nu+1}{K}\right)^{2\nu+1}}.\end{equation}
It is important to mention here that by using the Bernoulli inequality $(1+x)^a\geq 1+ax,$ for $x=-\frac{1}{K}>-1$ and $a=\nu+1>0,$ then clearly we have
$$\sqrt[4(\nu+1)]{\left(1-\frac{\nu+1}{K}\right)^{2\nu+1}}<\sqrt[4]{\left(1-\frac{1}{K}\right)^{2\nu+1}},$$
or in other words the inequality \eqref{ineq8} improves \eqref{ineq6}. To prove \eqref{ineq8} we would need to show that
\begin{equation}\label{ineqpsigen3}\Psi_{\nu}^{\frac{4(\nu+1)}{2\nu+1}}(x)+\frac{2(\nu+1)}{x}\Psi_{\nu}(x)-1<0.\end{equation}
By applying the identity \ref{quotgen} and the recurrence relation \ref{recgen} we obtain $$1-\frac{2(\nu+1)}{x}\Psi_{\nu}(x)=
\left(1+\frac{2(\nu+1)}{x}\frac{I_{\nu+1}(x)}{I_{\nu+2}(x)}\right)^{-1}=\frac{I_{\nu+2}(x)}{I_{\nu}(x)},$$
which implies that \eqref{ineqpsigen3} is equivalent to the Tur\'an type inequality
\begin{equation}\label{ineq9}I_{\nu+1}^{4\nu+4}(x)<I_{\nu+2}^{2\nu+1}(x)I_{\nu}^{2\nu+3}(x),\end{equation} which can be written as
$$\Psi_{\nu}^{2\nu+{3}}(x)<\Psi_{\nu+1}^{2\nu+{1}}(x),$$
where $x>0$ and $\nu\geq0.$ However, we were able to show the Tur\'an type inequality \eqref{ineq9} only for small values of $x.$ All the same, we believe that
\eqref{ineq9} is true for all $x>0$ and $\nu\geq0,$ and this open problem may of interest for further research. Using the recurrence relation
$$\Psi_{\nu}(x)\left(\frac{2(\nu+1)}{x}+\Psi_{\nu+1}(x)\right)=1$$
and the Amos bound $\Psi_{\nu}(x)>\Gamma_{\nu}(x),$ where $\nu\geq0,$ $x>0$ and \cite{amos}
$$\Gamma_{\nu}(x)=\frac{x}{\sqrt{x^2+\left(\nu+\frac{3}{2}\right)^{2}}+\nu+\frac{1}{2}},$$
we obtain that
$$\frac{\Psi_{\nu+1}^{2\nu+1}(x)}{\Psi_{\nu}^{2\nu+3}(x)}=\Psi_{\nu+1}^{2\nu+1}(x)\left(\frac{2(\nu+1)}{x}+\Psi_{\nu+1}(x)\right)^{2\nu+3}>
\Gamma_{\nu+1}^{2\nu+1}(x)\left(\frac{2(\nu+1)}{x}+\Gamma_{\nu+1}(x)\right)^{2\nu+3}>1$$
for $x\in(0,x_{\nu})$ and $\nu\geq0,$ where $x_{\nu}$ is the unique positive root of the equation
$$\Gamma_{\nu+1}^{2\nu+1}(x)\left(\frac{2(\nu+1)}{x}+\Gamma_{\nu+1}(x)\right)^{2\nu+3}=1.$$

Finally, we mention that the inequality \eqref{ineq8} is also sharp. Namely, if we consider another extension of $\lambda,$ that is,
\[\xi_\nu(x)=\frac{\log\left(1-\frac{2(\nu+1)}{x}\Psi_\nu(x)\right)}{\log\Psi_\nu(x)},\]
then the best constant $\gamma_\nu$ for which $\xi_{\nu}(x)<\gamma_{\nu}$ for $x>0$ and $\nu\geq0,$ equals to
\[\gamma_\nu=\lim_{x\to\infty}\xi_\nu(x)=\frac{4(\nu+1)}{2\nu+1}.\]
In particular, $\gamma_0=\beta=4.$ The steps here are also the same as in the special case above when $\nu=0.$ Recall that the asymptotic estimation \eqref{assym} yields that as $x$ grows we have \eqref{assymp2} and substituting this into the definition of $\xi_\nu(x)$ we can see that asymptotically it equals to
\[\xi_\nu(x)=\frac{\log(x-2(\nu+1)+O\left(\frac{1}{x}\right))-\log(x)}
{\log\left(8x-4(\nu+1)^2+1+O\left(\frac{1}{x}\right)\right)-\log\left(8x-4\nu^2+1+O\left(\frac{1}{x}\right)\right)},\]
and consequently
\[\xi_\nu(x)=\frac{\frac{2(\nu+1)}{x}+O\left(\frac{1}{x^2}\right)}{\left(\frac{(\nu+1)^2}{2}-\frac{\nu^2}{2}\right)\frac1x+O\left(\frac{1}{x^2}\right)}=
\frac{4(\nu+1)+O\left(\frac{1}{x}\right)}{2\nu+1+O\left(\frac{1}{x}\right)}\]
for large $x$, from which the result follows.

\section{Approximation of $r(K)$ from Lagrange's inversion and a rational approximation}
\setcounter{equation}{0}

In this section our aim is to propose two other approximations for the asymptotic order parameter of the stochastic Kuramoto model. The first one is
deduced by using the Lagrange inversion theorem, and the second one is a rational approximation.

\begin{center}
\begin{figure}[!ht]
   \centering
       \includegraphics[width=13cm]{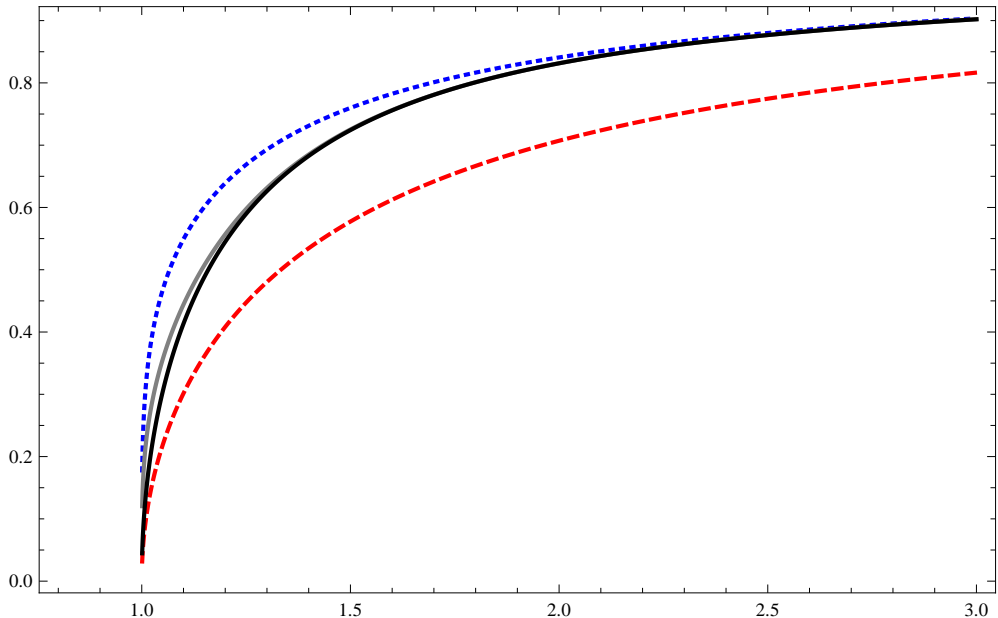}
       \caption{The graph of the function $K\mapsto r(K)$ on $[0,3]$ together with the bounds and the approximation.}
       \label{fig2}
\end{figure}
\end{center}

\subsection{Approximation from Lagrange inversion} We introduce the function $f:(0,\infty)\to\mathbb{R},$ defined by $f(r)=r-\Psi(2Kr).$ If $0<K\le1$, then there is only one non-negative real root of the equation $f(r)=0$ and this is $r=0$. But if $K>1$, there is an additional, non-trivial solution for this transcendental equation which we denote by $r(K)$. Recall the inequality \eqref{ineq1} and the fact that a sharper upper estimation is valid, see \eqref{ineq2}
\[r(K)<\sqrt[4]{1-\frac{1}{K}}.\]
Thanks to the analyticity of $f$ in the neighbor of its root, one can use the Lagrange inversion theorem in its simplest form to establish a better approximation to $r(K)$. This approximation for $K>1$ reads as
\begin{equation}
L(K)=A(K)+\frac{\Psi(s)-A(K)}{1-\Psi(s)/A(K)+2K\Psi^2(s)-2KI_2(s)/I_0(s)},\label{liestim}
\end{equation}
where
\[A(K)=\sqrt[4]{1-\frac{1}{K}},\quad\mbox{and}\quad s=2KA(K).\]
The above expression is nothing else but the zeroth order approximation plus the first order term in the inverse series of $f$ around the root estimation $A(K)$.
The performance of this approximation is drawn on Fig. \ref{fig2}. Here the upper, blue plot belongs to the upper bound $A(K)=\sqrt[4]{1-\frac1K}$, the bottom red plot belongs to the lower bound $\sqrt{1-\frac{1}{K}}$, the black one is the theoretical solution $r(K)$, while the grey plot is our approximation $L(K)$. One can see that for all $K>1$ the Lagrange estimation \eqref{liestim} approximates the theoretical $r(K)$ best. Numerical calculations show that if $K\ge2.8$ then $L(K)$ already gives 6 digits accuracy. That $L(K)$ is really better than $A(K)$ can easily be seen independently from the above graph. Indeed, the series defined by the Lagrange inversion theorem converges to the solution if the center is close enough to the theoretical solution. Since we started the approximation from the point $A(K)$, this latter requirement satisfies, and $A(K)$ is the zeroth order approximation. Then one more term in the Lagrange formula (resulting $L(K)$) gets even closer to $r(K)$.

\begin{center}
\begin{figure}[!ht]
   \centering
       \includegraphics[width=13cm]{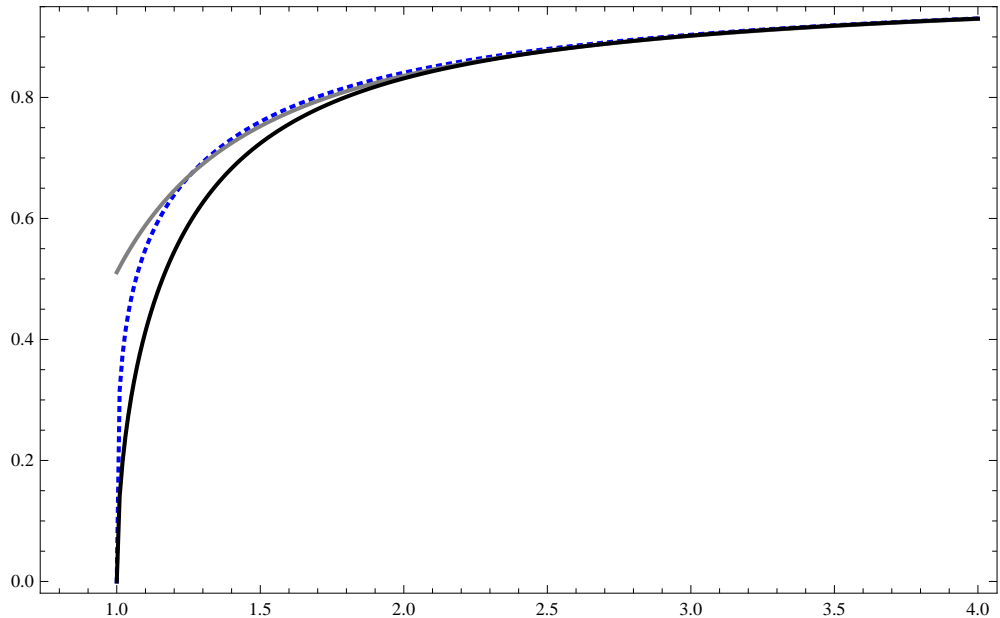}
       \caption{The graph of the function $K\mapsto r(K)$ on $[0,4]$ together with the rational approximation.}
       \label{fig3}
\end{figure}
\end{center}

\subsection{A rational approximation} The advantage of $A(K)$ is that it is algebraic. While $L(K)$ approximates the theoretical solution better, it is transcendental since it contains transcendental functions. In this section we find an approximation which is not simply algebraic but \emph{rational}.
We use the well known expansion \eqref{assym}, that is,
\[I_\nu(x)=\frac{e^x}{\sqrt{2\pi x}}\left(1-\frac{4\nu^2-1}{8x}+\frac{\left(4 \nu^2-1\right) \left(4 \nu^2-9\right)}{128x^2}-\frac{\left(4 \nu^2-1\right) \left(4 \nu^2-9\right) \left(4 \nu^2-25\right)}{1536x^3}+O\left(\frac{1}{x^4}\right)\right).\]
We truncate this at the ``$O$'' term, and substitute it into \eqref{liestim}. Thus, we get a fraction in which the nominator and denominator are polynomials of $K$ and $A(K)$. Since $A(K)\approx1$ as $K$ grows, we simply write 1 in place of $A(K)$. After a simplification we get the following expression
\[L_{\textrm{pol}}(K)=\frac{4 K \left(1048576 K^5-393216 K^4-276480 K^3+40320 K^2-7560 K-1575\right)}{4194304 K^6-524288 K^5-843776 K^4+376320 K^3+3936 K^2+540 K+3375}.\]

This is not, of course, better than $L(K)$. But $L_{\textrm{pol}}(K)$ offers a rather good rational approximation to $r(K)$ comparable at least with $A(K)$. The performance of $L_{\textrm{pol}}(K)$ is shown in Fig. \ref{fig3}. Here the blue plot represents $A(K)$, the grey line is of $L_{\textrm{pol}}(K)$, while the bottom black plot is of $r(K)$.
For some values we calculated the differences of $L_{\textrm{pol}}(K)$ and $A(K)$ with the theoretical solution

\begin{center}
\begin{tabular}{c|c|c|c|c|c}
$K$&1.5&2.0&5.0&10&100\\
$A(K)-r(K)$&0.035677&0.009434&0.0001994&0.00001936&$1.59\cdot10^{-8}$\\
$L_{\textrm{pol}}(K)-r(K)$&0.02818&0.0042565&$-0.000234$&$-0.0000372$&$-4.25\cdot10^{-8}$
\end{tabular}
\end{center}

\end{document}